\documentclass [12pt]{article}
\usepackage{amsmath,amssymb}
\begin{document}

\title{Analogs of Cramer's rule for the least squares
solutions of some matrix equations. }
\author{Ivan Kyrchei  \footnote{Pidstrygach Institute
for Applied Problems of Mechanics and Mathematics, str.Naukova 3b,
Lviv, 79005, Ukraine,  kyrchei@lms.lviv.ua}}
\date{}
 \maketitle


\begin{abstract} The least squares solutions  with the minimum norm
of the matrix equations ${\rm {\bf A}}{\rm {\bf X}} = {\rm {\bf
B}}$, ${\rm {\bf X}}{\rm {\bf A}} = {\rm {\bf B}}$ and ${\rm {\bf
A}}{\rm {\bf X}}{\rm {\bf B}} ={\rm {\bf D}} $ are considered in
this paper. We use the  determinantal representations of the Moore
- Penrose inverse obtained earlier by the author and  get analogs
of the Cramer rule for the least squares solutions of these matrix
equations.
\end{abstract}

\textit{Keywords}: Moore-Penrose inverse,   matrix equation, least
squares solution, Cramer rule

\noindent \textit{MSC 2010}: 15A09, 15A24

\section{Introduction}
\newtheorem{definition}{Definition}[section]
\newtheorem{lemma}{Lemma}[section]
 \newtheorem{corollary}{Corollary}[section]
\newtheorem{theorem}{Theorem}[section]
\newtheorem{remark}{Remark}[section]

In this paper we shall adopt the following notation. Let ${\mathbb
C}^{m\times n} $ be the set of $m$ by $n$ matrices with complex
entries, ${\mathbb C}^{m\times n}_{r} $ be a subset of ${\mathbb
C}^{m\times n} $  in which any matrix has rank $r$, ${\rm \bf
I}_{m}$ be the identity matrix of order $m$, and $\|.\|$ be the
Frobenius norm of a matrix.

Denote by ${\rm {\bf a}}_{.j} $ and ${\rm {\bf a}}_{i.} $ the
$j$th column  and the $i$th row of ${\rm {\bf A}}\in {\mathbb
C}^{m\times n} $, respectively. Then   ${\rm {\bf a}}^{\ast}_{.j}
$ and ${\rm {\bf a}}^{\ast}_{i.} $ denote the $j$th column  and
the $i$th row of a Hermitian adjoint matrix ${\rm {\bf A}}^{\ast}$
as well. Let ${\rm {\bf A}}_{.j} \left( {{\rm {\bf b}}} \right)$
denote the matrix obtained from ${\rm {\bf A}}$ by replacing its
$j$th column with the  vector ${\rm {\bf b}}$, and by ${\rm {\bf
A}}_{i.} \left( {{\rm {\bf b}}} \right)$ denote the matrix
obtained from ${\rm {\bf A}}$ by replacing its $i$th row with
${\rm {\bf b}}$.

Let $\alpha : = \left\{ {\alpha _{1} ,\ldots ,\alpha _{k}}
\right\} \subseteq {\left\{ {1,\ldots ,m} \right\}}$ and $\beta :
= \left\{ {\beta _{1} ,\ldots ,\beta _{k}}  \right\} \subseteq
{\left\{ {1,\ldots ,n} \right\}}$ be subsets of the order $1 \le k
\le \min {\left\{ {m,n} \right\}}$. Then ${\left| {{\rm {\bf
A}}_{\beta} ^{\alpha} } \right|}$ denotes the minor of ${\rm {\bf
A}}$ determined by the rows indexed by $\alpha$ and the columns
indexed by $\beta$. Clearly, ${\left| {{\rm {\bf A}}_{\alpha}
^{\alpha} } \right|}$ be a principal minor determined by the rows
and  columns indexed by $\alpha$. For $1 \leq k\leq n$, denote by
 \[\textsl{L}_{ k, n}: = {\left\{
{\,\alpha :\alpha = \left( {\alpha _{1} ,\ldots ,\alpha _{k}}
\right),\,{\kern 1pt} 1 \le \alpha _{1} \le \ldots \le \alpha _{k}
\le n} \right\}}\]
 the collection of strictly increasing sequences
of $k$ integers chosen from the set $\left\{ {1,\ldots ,n}
\right\}$. For fixed $i \in \alpha $ and $j \in \beta $, let
\[I_{r,\,m} {\left\{ {i} \right\}}: = {\left\{ {\,\alpha :\alpha
\in L_{r,m} ,i \in \alpha} \right\}}{\rm ,} \quad J_{r,\,n}
{\left\{ {j} \right\}}: = {\left\{ {\,\beta :\beta \in L_{r,n} ,j
\in \beta} \right\}}.\]

Matrix equation is one of the important study fields of linear
algebra. Linear matrix equations, such as
\begin{equation}\label{eq:AX=B}
 {\rm {\bf A}}{\rm {\bf X}} = {\rm {\bf B}},
\end{equation}
\begin{equation}\label{eq:XB=D}
 {\rm {\bf X}}{\rm {\bf B}} = {\rm {\bf D}},
\end{equation}
 and \begin{equation}\label{eq:AXB=D}
 {\rm {\bf A}}{\rm {\bf X}}{\rm {\bf B}} = {\rm {\bf
D}},
\end{equation}
play an important role in linear system theory therefore a large
number of papers have presented several methods for solving these
matrix equations \cite{da,de,do,pe,vet}. In \cite{kh}, Khatri and
Mitra studied the Hermitian solutions to the  matrix equations
(\ref{eq:AX=B}) and (\ref{eq:AXB=D}) over the complex field and
the system of the equations  (\ref{eq:AX=B}) and (\ref{eq:XB=D}).
Wang, in \cite{wa1,wa2}, and Li and Wu, in \cite{li}   studied the
bisymmetric, symmetric and skew-antisymmetric least squares
solution to this system over the  quaternion skew field. Extreme
ranks of real matrices in least squares solution of the equation
(\ref{eq:AXB=D}) was investigated in \cite{liu} over the complex
field and in \cite{wa3} over the  quaternion skew field.

As we know, the Cramer rule gives an explicit expression for the
solution of nonsingular linear equations. In \cite{ro}, Robinson
gave its elegant proof  over the complex  field which aroused
great interest in finding determinantal formulas as analogs of the
Cramer rule for the matrix equations
\cite{be,ch,ji,ky,ky1,ky2,ve,wagb,wei,wer}.
 The Cramer rule for solutions of the
restricted matrix equations (\ref{eq:AX=B}), (\ref{eq:XB=D}) and
(\ref{eq:AXB=D}) was established in \cite{wag}.

In this paper, we use the results of \cite{ky} to obtain the
Cramer rule for least squares solutions of the matrix equations
(\ref{eq:AX=B}), (\ref{eq:XB=D}) and (\ref{eq:AXB=D}) without any
restriction. The paper is organized as follows. We start with some
basic concepts and results about determinantal representations
 of the Moore - Penrose inverse in Section 2. In Section 3, we derive some generalized Cramer rules for the
matrix equations (\ref{eq:AX=B}), (\ref{eq:XB=D}) and
(\ref{eq:AXB=D}). In Section 4, we show a numerical example to
illustrate the main result.
\section{Determinantal representations
 of the Moore - Penrose inverse}
\begin{definition}
If ${\rm {\bf A}}\in {\mathbb C}^{m\times n} $, then the matrix
${\rm {\bf A}}^{ +} $ is called the Moore - Penrose  inverse of
${\rm {\bf A}}$ if it satisfies the equations:
  1)\,$ \left( {{\rm {\bf A}}{\rm {\bf
A}}^{ +} } \right)^{ *}  = {\rm {\bf A}}{\rm {\bf A}}^{ +};$
  2)\,$ \left( {{\rm {\bf A}}^{ +} {\rm {\bf A}}} \right)^{ *}  = {\rm
{\bf A}}^{ +} {\rm {\bf A}};$
  3)\,$ {\rm {\bf A}}{\rm {\bf A}}^{ +}
{\rm {\bf A}} = {\rm {\bf A}};$
  4)\,${\rm {\bf A}}^{ +} {\rm {\bf
A}}{\rm {\bf A}}^{ +}  = {\rm {\bf A}}^{ +}. $
\end{definition}
It is well known the following proposition.
\begin{lemma}  For an arbitrary
${\rm {\bf A}}\in {\mathbb C}^{m\times n} $ there exists a unique
Moore - Penrose inverse ${\rm {\bf A}}^{ + }$.
\end{lemma}
\begin{lemma}\cite{ca}\label{kyrc2}
If ${\rm {\bf A}}\in {\mathbb C}^{m\times n} $, then
\[{\rm {\bf
A}}^{ +} = {\mathop {\lim }\limits_{\lambda \to 0}} {\rm {\bf
A}}^{ *} \left( {{\rm {\bf A}}{\rm {\bf A}}^{ *}  + \lambda {\rm
{\bf I}}} \right)^{ - 1} = {\mathop {\lim }\limits_{\lambda\to 0}}
\left( {{\rm {\bf A}}^{ *} {\rm {\bf A}} + \lambda {\rm {\bf I}}}
\right)^{ - 1}{\rm {\bf A}}^{ *}, \]
 where $\lambda \in {\mathbb
R} _{ +}  $, and ${\mathbb R} _{ +} $ is a set of the real
positive numbers.
\end{lemma}

\begin{corollary}\label{cor:A+A-1} If
${\rm {\bf A}}\in {\mathbb C}^{m\times n} $, then the following
statements are true.
 \begin{itemize}
\item [ i)] If $\rm{rank}\,{\rm {\bf A}} = n$, then ${\rm {\bf A}}^{ +}
= \left( {{\rm {\bf A}}^{ *} {\rm {\bf A}}} \right)^{ - 1}{\rm
{\bf A}}^{ * }$ .
\item [ ii)] If $\rm{rank}\,{\rm {\bf A}} =
m$, then ${\rm {\bf A}}^{ +}  = {\rm {\bf A}}^{ * }\left( {{\rm
{\bf A}}{\rm {\bf A}}^{ *} } \right)^{ - 1}.$
\item [ iii)] If $\rm{rank}\,{\rm {\bf A}} = n = m$, then ${\rm {\bf
A}}^{ +}  = {\rm {\bf A}}^{ - 1}$ .
\end{itemize}
\end{corollary}
\begin{theorem}\cite{ky}\label{theor:det_repr_MP}
If ${\rm {\bf A}} \in {\rm {\mathbb{C}}}_{r}^{m\times n} $ and
$r<min\{m,n\}$, then the Moore-Penrose inverse  ${\rm {\bf A}}^{
+} = \left( {a_{ij}^{ +} } \right) \in {\rm
{\mathbb{C}}}_{}^{n\times m} $ possess the following determinantal
representations:
\begin{equation}
\label{eq:det_repr_A*A}
 a_{ij}^{ +}  = {\frac{{{\sum\limits_{\beta
\in J_{r,\,n} {\left\{ {i} \right\}}} {{\left| \left( {\left(
{{\rm {\bf A}}^{ *} {\rm {\bf A}}} \right)_{\,. \,i} \left( {{\rm
{\bf a}}_{.j}^{ *} }  \right)} \right){\kern 1pt} {\kern 1pt}
_{\beta} ^{\beta} \right|}} }  }}{{{\sum\limits_{\beta \in
J_{r,\,\,n}} {{\left| {\left( {{\rm {\bf A}}^{ *} {\rm {\bf A}}}
\right){\kern 1pt} _{\beta} ^{\beta} }  \right|}}} }}},
\end{equation}
or
\begin{equation}
\label{eq:det_repr_AA*} a_{ij}^{ +}  =
{\frac{{{\sum\limits_{\alpha \in I_{r,m} {\left\{ {j} \right\}}}
{{\left| \left( {({\rm {\bf A}}{\rm {\bf A}}^{ *} )_{j\,.\,} ({\rm
{\bf a}}_{i.\,}^{ *} )} \right)\,_{\alpha} ^{\alpha}\right|} }
}}}{{{\sum\limits_{\alpha \in I_{r,\,m}}  {{\left| {\left( {{\rm
{\bf A}}{\rm {\bf A}}^{ *} } \right){\kern 1pt} _{\alpha}
^{\alpha} } \right|}}} }}}.
\end{equation}
\end{theorem}
\begin{remark}\label{rem:det_repr_MP_A*A}
If ${\rm rank}\,{\rm {\bf A}} = n$, then by Corollary
\ref{cor:A+A-1} ${\rm {\bf A}}^{ +}  = \left( {{\rm {\bf A}}^{ *}
{\rm {\bf A}}} \right)^{ - 1}{\rm {\bf A}}^{ *} $. Therefore, we
get the following representation of  ${\rm {\bf A}}^{ +} $:

\begin{equation}
\label{eq:det_repr_MP_A*A}{\rm {\bf A}}^{ +}  = {\frac{{1}}{\det
({{\rm {\bf A}}^{ *} {\rm {\bf A}}})}}
\begin{pmatrix}
  {\det ({{\rm {\bf A}}^{ *} {\rm {\bf
A}}})_{.\,1} \left( {{\rm {\bf a}}_{.\,1}^{ *} }  \right)} &
\ldots & {\det ({{\rm {\bf A}}^{ *} {\rm {\bf A}}})_{.\,1} \left(
{{\rm {\bf a}}_{.\,m}^{ *} }  \right)} \\
  \ldots & \ldots & \ldots \\
  {\det ({{\rm {\bf A}}^{ *} {\rm {\bf A}}})_{.\,n}
 \left( {{\rm {\bf a}}_{.\,1}^{ *} }  \right)} & \ldots &
 {\det ({{\rm {\bf A}}^{ *} {\rm {\bf
A}}})_{.\,n} \left( {{\rm {\bf a}}_{.\,m}^{ *} }  \right)}.
\end{pmatrix}
\end{equation}
If ${\rm rank}\,{\rm {\bf A}} = n<m$, then by Theorem
\ref{theor:det_repr_MP} for ${\rm {\bf A}}^{ +} $ we have
(\ref{eq:det_repr_A*A}) as well.
\end{remark}
\begin{remark}\label{rem:det_repr_MP_AA*}
If ${\rm rank}\,{\rm {\bf A}} = m$, then by Corollary
\ref{cor:A+A-1} ${\rm {\bf A}}^{ +}  = {\rm {\bf A}}^{ *} \left(
{{\rm {\bf A}}{\rm {\bf A}}^{ *} } \right)^{ - 1}$. Hence, we
obtain the following representation of ${\rm {\bf A}}^{ +} $:

\begin{equation}
\label{eq:det_repr_MP_AA*} {\rm {\bf A}}^{ +}  = {\frac{{1}}{{\det
({\rm {\bf A}}{\rm {\bf A}}^{ *} )}}}\begin{pmatrix}
  {\det ({\rm {\bf A}}{\rm {\bf A}}^{ *} )_{1\,.} \left( {{\rm {\bf a}}_{1\,.}^{ *} }  \right)}
   & \ldots &
   {\det ({\rm {\bf A}}{\rm {\bf A}}^{ *} )_{m\,.} \left( {{\rm {\bf a}}_{1\,.}^{ *} }  \right)} \\
 \ldots & \ldots & \ldots \\
 {\det ({\rm {\bf A}}{\rm {\bf A}}^{ *} )_{1\,.} \left( {{\rm {\bf a}}_{n\,.}^{ *} }  \right)}
  & \ldots &
{\det ({\rm {\bf A}}{\rm {\bf A}}^{ *} )_{m\,.} \left( {{\rm {\bf
a}}_{n\,.}^{ *} }  \right)}
\end{pmatrix}
\end{equation}
If ${\rm rank}\,{\rm {\bf A}} = m < n$, then by Theorem
\ref{theor:det_repr_MP} for ${\rm {\bf A}}^{ +} $ we also have
(\ref{eq:det_repr_AA*}).
\end{remark}
\section{Cramer's rule of the least squares solution of some  matrix
equation}
\begin{definition}
Consider a  matrix equation (\ref{eq:AX=B}), where ${\rm {\bf
A}}\in {\mathbb{C}}^{m\times n},{\rm {\bf B}} \in
{\mathbb{C}}^{m\times s} $ are given, ${\rm {\bf X}} \in
{\mathbb{C}}^{n\times s}$ is unknown. Suppose
\[S_{1}=\{{\rm {\bf X}}|{\rm {\bf X}} \in
{\mathbb{C}}^{n\times s}, \|{\rm {\bf A}}{\rm {\bf X}} - {\rm {\bf
B}}\|=\min \}.\] Then matrices ${\rm {\bf X}} \in
{\mathbb{C}}^{n\times s}$ such that ${\rm {\bf X}}\in S_{1}$ are
called least squares solutions of the matrix equation
(\ref{eq:AX=B}).
 If
${\rm {\bf X}}_{LS}={\min}_ {{\rm {\bf X}}\in S_{1}}\|{\rm {\bf
X}}\|$, then ${\rm {\bf X}}_{LS}$ is called the least squares
solution of (\ref{eq:AX=B}) with minimum norm.
\end{definition}
If  the equation (\ref{eq:AX=B}) has no precision solutions, then
${\rm {\bf X}}_{LS}$ is its maximum approximate solution.

The following important theorem is well-known.

\begin{theorem}(\cite{ji})\label{theor:LS} The least squares solutions of (\ref{eq:AX=B}) are
\[{\rm {\bf X}} = {\rm {\bf A}}^{+}{\rm {\bf B}} + ({\rm {\bf
I}}_{n} - {\rm {\bf A}}^{+} {\rm {\bf A}}){\rm {\bf C}},\] in
which ${\rm {\bf C}}\in {\mathbb{C}}^{n\times s}$ is an arbitrary
 matrix and the minimum norm solution is ${\rm {\bf
X}}_{LS}={\rm {\bf A}}^{+}{\rm {\bf B}}$.
\end{theorem}
 We denote ${\rm {\bf A}}^{ \ast}{\rm {\bf B}}=:\hat{{\rm
{\bf B}}}= (\hat{b}_{ij})\in {\mathbb{C}}^{n\times s}$.
\begin{theorem}
\begin{enumerate}
\item[(i)] If ${\rm rank}\,{\rm {\bf A}} = n$, then for the
least squares solution ${\rm {\bf X}}_{LS}=(x_{ij})\in
{\mathbb{C}}^{n\times s}$ for all $i = \overline {1,n} $, $j =
\overline {1,s} $ we have
\begin{equation}
\label{eq:ls_AX_full} x_{i\,j} = {\frac{{\det ({\rm {\bf A}}^{ *}
{\rm {\bf A}})_{.\,i} \left( {\hat{{\rm {\bf b}}}_{.j}}
\right)}}{{ \det ({\rm {\bf A}}^{ *} {\rm {\bf A}})}}}
\end{equation}
\noindent where $\hat{{\rm {\bf b}}}_{.j}$ is the $j$th column of
$\hat{{\rm {\bf B}}}$ for all $j = \overline {1,s} $.
\item[(ii)] If ${\rm rank}\,{\rm {\bf A}} = r \le m < n$, then for all
$i = \overline {1,n} $, $j = \overline {1,s} $ we have
\begin{equation}
\label{eq:ls_AX} x_{ij} = {\frac{{{\sum\limits_{\beta \in
J_{r,\,n} {\left\{ {j} \right\}}} {{\left| \left( {\left( {{\rm
{\bf A}}^{ *} {\rm {\bf A}}} \right)_{\,.\,i} \left( {\hat{{\rm
{\bf b}}}_{.j}} \right)} \right){\kern 1pt} {\kern 1pt} _{\beta}
^{\beta} \right|}} } }}{{{\sum\limits_{\beta \in J_{r,\,\,n}}
{{\left| {\left( {{\rm {\bf A}}^{ *} {\rm {\bf A}}} \right){\kern
1pt} {\kern 1pt} _{\beta} ^{\beta} }  \right|}}} }}}.
\end{equation}
\end{enumerate}
\end{theorem}
{\textit{Proof.}} i) If ${\rm rank}\,{\rm {\bf A}} = n$,  then by
Corollary \ref{cor:A+A-1} ${\rm {\bf A}}^{ +}  = \left( {{\rm {\bf
A}}^{ *} {\rm {\bf A}}} \right)^{ - 1}{\rm {\bf A}}^{*} $.  By
multiplying $ {\rm {\bf A}}^{+}{\rm {\bf B}}=\left( {{\rm {\bf
A}}^{*} {\rm {\bf A}}} \right)^{ - 1}{\rm {\bf A}}^{*}{\rm {\bf
B}}=\left( {{\rm {\bf A}}^{*} {\rm {\bf A}}} \right)^{-
1}{\hat{{\rm {\bf B}}}} $, we obtain for all $i = \overline {1,n}
$, $j = \overline {1,s} $
\[
x_{ij}  = {\frac{{1}}{{\det ({\rm {\bf A}}^{ *} {\rm {\bf
A}})}}}{\sum\limits_{k = 1}^{n} {{L}_{ki} \hat{b}_{kj}} },
\]
\noindent where ${L}_{ij}$ is a  $ij$th
 cofactor of $({\rm {\bf A}}^{ \ast}{\rm {\bf A}})$ for all $i,j = \overline {1,n}
$. Denoting the $j$th column of $\hat{{\rm {\bf B}}}$ by
$\hat{{\rm {\bf b}}}_{.j} $, it follows (\ref{eq:ls_AX_full}).

\noindent ii) If ${\rm rank}\,{\rm {\bf A}} = r \le m < n$, then
by Theorem  \ref{theor:det_repr_MP} we can represent the matrix
${\rm {\bf A}}^{ +} $ by (\ref{eq:det_repr_A*A}). Therefore, we
obtain for all $i = \overline {1,n} $, $j = \overline {1,s} $
\[
x_{ij}  =\sum_{k=1}^{m}
a_{ik}^{+}b_{kj}=\sum_{k=1}^{m}{\frac{{{\sum\limits_{\beta \in
J_{r,\,n} {\left\{ {i} \right\}}} {{\left| \left( {\left( {{\rm
{\bf A}}^{ *} {\rm {\bf A}}} \right)_{\,. \,i} \left( {{\rm {\bf
a}}_{.k}^{ *} }  \right)} \right) {\kern 1pt} _{\beta}
^{\beta}\right|} } } }}{{{\sum\limits_{\beta \in J_{r,\,\,n}}
{{\left| {\left( {{\rm {\bf A}}^{ *} {\rm {\bf A}}} \right){\kern
1pt} _{\beta} ^{\beta} }  \right|}}} }}}\cdot b_{kj}=
\]
\[
{\frac{{{\sum\limits_{\beta \in J_{r,\,n} {\left\{ {i}
\right\}}}\sum_{k=1}^{m} {{\left| \left( {\left( {{\rm {\bf A}}^{
*} {\rm {\bf A}}} \right)_{\,. \,i} \left( {{\rm {\bf a}}_{.k}^{
*} } \right)} \right) {\kern 1pt} _{\beta} ^{\beta}\right|} } }
}\cdot b_{kj}} {{{\sum\limits_{\beta \in J_{r,\,\,n}} {{\left|
{\left( {{\rm {\bf A}}^{ *} {\rm {\bf A}}} \right){\kern 1pt}
_{\beta} ^{\beta} } \right|}}} }}}.
\]
Since ${\sum\limits_{k} {{\rm {\bf a}}_{.\,k}^{ *}  b_{kj}} }=
\left( {{\begin{array}{*{20}c}
 {{\sum\limits_{k} {a_{1k}^{ *}  b_{kj}} } } \hfill \\
 {{\sum\limits_{k} {a_{2k}^{ *}  b_{kj}} } } \hfill \\
 { \vdots}  \hfill \\
 {{\sum\limits_{k} {a_{nk}^{ *}  b_{kj}} } } \hfill \\
\end{array}} } \right) = \hat{{\rm {\bf b}}}_{.j}$, then it follows (\ref{eq:ls_AX}).
$\blacksquare$
\begin{definition}
Consider a  matrix equation
\begin{equation}\label{eq:XA=B}
 {\rm {\bf X}}{\rm {\bf A}} = {\rm {\bf B}},
\end{equation}
\noindent  where ${\rm {\bf A}}\in {\mathbb{C}}^{m\times n},{\rm
{\bf B}} \in {\mathbb{C}}^{s\times n} $ are given, ${\rm {\bf X}}
\in {\mathbb{C}}^{s\times m}$ is unknown. Suppose
\[S_{2}=\{{\rm {\bf X}}|\,{\rm {\bf X}} \in
{\mathbb{C}}^{s\times m}, \|{\rm {\bf X}}{\rm {\bf A}} - {\rm {\bf
B}}\|=\min \}.\] Then matrices ${\rm {\bf X}} \in
{\mathbb{C}}^{s\times m}$ such that ${\rm {\bf X}}\in S_{2}$ are
called least squares solutions of the matrix equation
(\ref{eq:XA=B}).
 If
${\rm {\bf X}}_{LS}={\min}_ {{\rm {\bf X}}\in S_{2}}\|{\rm {\bf
X}}\|$, then ${\rm {\bf X}}_{LS}$ is called the least squares
solution of (\ref{eq:XA=B}) with minimum norm.
\end{definition}
The following  theorem can be obtained by analogy to Theorem
\ref{theor:LS}.

\begin{theorem} The least squares solutions of (\ref{eq:XA=B}) are
\[{\rm {\bf X}} = {\rm {\bf B}}{\rm {\bf A}}^{+} + {\rm {\bf C}}({\rm {\bf
I}}_{m} -  {\rm {\bf A}}{\rm {\bf A}}^{+}),\] in which ${\rm {\bf
C}}\in {\mathbb{C}}^{s\times m}$ is an arbitrary matrix and the
minimum norm solution is ${\rm {\bf X}}_{LS}={\rm {\bf B}}{\rm
{\bf A}}^{+}$.
\end{theorem}
 We denote  ${\rm {\bf B}}{\rm {\bf A}}^{
\ast}=:\check{{\rm {\bf B}}}= (\check{b}_{ij})\in
{\mathbb{C}}^{s\times m}$.
\begin{theorem}
\begin{enumerate}
\item[(i)] If ${\rm rank}\,{\rm {\bf A}} = m$, then for the
least squares solution ${\rm {\bf X}}_{LS}=(x_{ij})\in
{\mathbb{C}}^{s\times m}$ for all $i = \overline {1,s} $, $j =
\overline {1,m} $ we have
\begin{equation}
\label{eq:ls_XA_full} x_{i\,j} = {\frac{{\det ({\rm {\bf A}}{\rm
{\bf A}}^{ *} )_{j.\,} \left( {\check{{\rm {\bf b}}}_{i\,.}}
\right)}}{{\det ({\rm {\bf A}}{\rm {\bf A}}^{ *} )}}}
\end{equation}
\noindent where $\check{{\rm {\bf b}}}_{i.}$ is the $i$th row of
$\check{{\rm {\bf B}}}$ for all $i = \overline {1,s} $.
\item[(ii)] If ${\rm rank}\,{\rm {\bf A}} = r \le n < m$, then for all $i = \overline {1,s} $, $j =
\overline {1,m} $ we have
\begin{equation}
\label{eq:ls_XA} x_{ij} = {\frac{{{\sum\limits_{\alpha \in I_{r,m}
{\left\{ {i} \right\}}} {{\left| \left( {\left( {{\rm {\bf A}}{\rm
{\bf A}}^{ *} } \right)_{\,j\,.} \left( {\check{{\rm {\bf
b}}}_{i\,.}} \right)} \right)\,_{\alpha} ^{\alpha}\right|} }
}}}{{{\sum\limits_{\alpha \in I_{r,\,m}}  {{\left| {\left( {{\rm
{\bf A}}{\rm {\bf A}}^{ *} } \right) {\kern 1pt} _{\alpha}
^{\alpha} } \right|}}} }}}.
\end{equation}
\end{enumerate}
\end{theorem}
{\textit{Proof.}} i) If ${\rm rank}\,{\rm {\bf A}} = m$,  then by
Corollary \ref{cor:A+A-1} ${\rm {\bf A}}^{ +}  = {\rm {\bf A}}^{ *
}\left( {{\rm {\bf A}}{\rm {\bf A}}^{ *} } \right)^{ - 1}.$   By
multiplying $ {\rm {\bf B}}{\rm {\bf A}}^{+}={\rm {\bf B}}{\rm
{\bf A}}^{ * }\left( {{\rm {\bf A}}{\rm {\bf A}}^{ *} } \right)^{
- 1}={\check{{\rm {\bf B}}}}\left( {{\rm {\bf A}}{\rm {\bf A}}^{
*} } \right)^{ - 1} $, we obtain for all $i = \overline {1,s} $,
$j = \overline {1,m} $,
\[
x_{ij}  = {\frac{{1}}{{\det ({\rm {\bf A}}{\rm {\bf A}}^{ *}
)}}}{\sum\limits_{k = 1}^{m}\check{b}_{ik} {{R}_{jk} } },
\]
\noindent where $R_{ij}$ is a  $ij$th
 cofactor of $({\rm {\bf A}}{\rm {\bf A}}^{ \ast})$ for all $i,j =
 \overline {1,m}$. Denoting the $i$th row of
$\check{{\rm {\bf B}}}$ by $\check{{\rm {\bf b}}}_{i.}$, it
follows (\ref{eq:ls_XA_full}).

\noindent ii) If ${\rm rank}\,{\rm {\bf A}} = r \le n <m $, then
by Theorem  \ref{theor:det_repr_MP} we can represent the matrix
${\rm {\bf A}}^{ +} $ by (\ref{eq:det_repr_AA*}). Therefore, for
all $i = \overline {1,s} $, $j = \overline {1,m} $ we obtain
\[
x_{ij}  =\sum_{k=1}^{n}b_{ik} a_{kj}^{+}=\sum_{k=1}^{n}b_{ik}\cdot
{\frac{{{\sum\limits_{\alpha \in I_{r,m} {\left\{ {i} \right\}}}
{{\left| \left( {\left( {{\rm {\bf A}}{\rm {\bf A}}^{ *} }
\right)_{\,j\,.} \left( {\bf a}_{k\,.}^{*} \right)}
\right)\,_{\alpha} ^{\alpha}\right|} } }}}{{{\sum\limits_{\alpha
\in I_{r,\,m}}  {{\left| {\left( {{\rm {\bf A}}{\rm {\bf A}}^{ *}
} \right) {\kern 1pt} _{\alpha} ^{\alpha} } \right|}}} }}} =
\]
\[
{\frac{{{\sum_{k=1}^{n}b_{ik} \sum\limits_{\alpha \in I_{r,m}
{\left\{ {i} \right\}}} {{\left| \left( {\left( {{\rm {\bf A}}{\rm
{\bf A}}^{ *} } \right)_{\,j\,.} \left( {\bf a}_{k\,.}^{*}
\right)} \right)\,_{\alpha} ^{\alpha}\right|} } } }}
{{{\sum\limits_{\alpha \in I_{r,\,m}}  {{\left| {\left( {{\rm {\bf
A}}{\rm {\bf A}}^{ *} } \right) {\kern 1pt} _{\alpha} ^{\alpha} }
\right|}}} }}}
\]
Since for all $i = \overline {1,s} $
 \[{\sum\limits_{k} {{  b_{ik}\rm {\bf a}}_{k\,.}^{ *}}
}=\begin{pmatrix}
  \sum\limits_{k} {b_{ik}a_{k1}^{ *} } & \sum\limits_{k} {b_{ik}a_{k2}^{ *} } & \cdots & \sum\limits_{k} {b_{ik}a_{km}^{ *} }
\end{pmatrix}
 = \check{{\rm {\bf b}}}_{i.},\]
then it follows (\ref{eq:ls_XA}). $\blacksquare$

\begin{definition}
Consider a  matrix equation (\ref{eq:AXB=D}),  where $ {\rm {\bf
A}}\in{\rm {\mathbb{C}}}^{m \times n}_{r_{1}},{\rm {\bf B}}\in{\rm
{\mathbb{C}}}^{p \times q}_{r_{2}}, {\rm {\bf D}}\in{\rm
{\mathbb{C}}}^{m \times q}$ are given, $ {\rm {\bf X}}\in{\rm
{\mathbb{C}}}^{n \times p}$ is unknown. Suppose
\[S_{3}=\{{\rm {\bf X}}|\,{\rm {\bf X}} \in
{\mathbb{C}}^{n \times p}, \|{\rm {\bf A}}{\rm {\bf X}}{\rm {\bf
B}} - {\rm {\bf D}}\|=\min \}.\] Then matrices ${\rm {\bf X}} \in
{\mathbb{C}}^{n \times p}$ such that ${\rm {\bf X}}\in S_{3}$ are
called least squares solutions of the matrix equation
(\ref{eq:AXB=D}).
 If
${\rm {\bf X}}_{LS}={\min}_ {{\rm {\bf X}}\in S_{3}}\|{\rm {\bf
X}}\|$, then ${\rm {\bf X}}_{LS}$ is called the least squares
solution of (\ref{eq:AXB=D}) with minimum norm.
\end{definition}
The following important theorem is well-known.
\begin{theorem}(\cite{wer})\label{theor:LS_AXB} The least squares solutions of (\ref{eq:AXB=D}) are
\[{\rm {\bf X}} = {\rm {\bf A}}^{+}{\rm {\bf
D}}{\rm {\bf B}}^{+} + ({\rm {\bf I}}_{n} - {\rm {\bf A}}^{+} {\rm
{\bf A}}){\rm {\bf V}}+{\rm {\bf W}}({\rm {\bf I}}_{p} -{\rm{\bf
B}}{\rm{\bf B}}^{+}),\] in which $\{{\rm {\bf V}},{\rm {\bf
W}}\}\subset {\mathbb{C}}^{n\times p}$  are  arbitrary quaternion
matrices and the least squares solution with minimum norm is ${\rm
{\bf X}}_{LS}={\rm {\bf A}}^{+}{\rm {\bf D}}{\rm {\bf B}}^{+}$.
\end{theorem}
 We denote  ${\rm {\bf \widetilde{D}}}= {\rm {\bf
A}}^\ast{\rm {\bf D}}{\rm {\bf B}}^\ast$.
\begin{theorem}\label{theor:AXB=D}
\begin{enumerate}
\item[(i)] If ${\rm rank}\,{\rm {\bf A}} = n$ and ${\rm rank}\,{\rm {\bf B}} = p$, then for the
least squares solution ${\rm {\bf X}}_{LS}=(x_{ij})\in
{\mathbb{C}}^{n\times p}$  of (\ref{eq:AXB=D}) we have for all $i
= \overline {1,n} $, $j = \overline {1,p} $,
\begin{equation}
\label{eq:AXB_cdetA*A} x_{i\,j} = {\frac{{\det \left( {({\rm {\bf
A}}^{ *} {\rm {\bf A}})_{.\,i\,} \left( {{\rm {\bf d}}_{.j}^{{\rm
{\bf B}}}} \right)}\right)}}{{\det ({\rm {\bf A}}^{ *} {\rm {\bf
A}})\cdot \det ({\rm {\bf B}}{\rm {\bf B}}^{ *} )}}},
\end{equation}
or
\begin{equation}
\label{eq:AXB_rdetBB*} x_{i\,j} = {\frac{{\det\left( { ({\rm {\bf
B}}{\rm {\bf B}}^{ *} )_{j.\,} \left( {{\rm {\bf d}}_{i\,.}^{{\rm
{\bf A}}}} \right)}\right)}}{{\det ({\rm {\bf A}}^{ *} {\rm {\bf
A}})\cdot \det ({\rm {\bf B}}{\rm {\bf B}}^{ *} )}}},
\end{equation}
 \noindent where
  \begin{equation}
\label{eq:def_d^B}{\rm {\bf d}}_{.j}^{{\rm {\bf B}}} : = \left[
{\det\left( ({\rm {\bf B}}{\rm {\bf B}}^{ *} )_{j.\,} \left(
{\tilde{{\rm {\bf d}}}_{1\,.}^{}} \right)\right),\ldots
,\det\left( ({\rm {\bf B}}{\rm {\bf B}}^{ *} )_{j.\,} \left(
{\tilde{{\rm {\bf d}}}_{n\,.}^{}} \right)\right)} \right]^{T},
\end{equation}
   \begin{equation}
\label{eq:def_d^A}{\rm {\bf d}}_{i\,.}^{{\rm {\bf A}}} : = \left[
{\det  \left(({\rm {\bf A}}^{ *} {\rm {\bf A}})_{.\,i\,} \left(
{\tilde{{\rm {\bf d}}}_{.1}} \right)\right),\ldots ,\det \left(
({\rm {\bf A}}^{ *} {\rm {\bf A}})_{.\,i\,} \left( {\tilde{{\rm
{\bf d}}}_{.p}} \right)\right)} \right]
\end{equation}
are respectively the column-vector and the row-vector.
$\tilde{{\rm {\bf d}}}_{i\,.}$ is the  $i$th row of ${\rm {\bf
\widetilde{D}}}$ for all $i = \overline {1,n} $, and $\tilde{{\rm
{\bf d}}}_{.\,j}$ is the $j$th column of ${\rm {\bf
\widetilde{D}}}$ for all $j = \overline {1,p} $.

\item[(ii)] If ${\rm rank}\,{\rm {\bf A}} = r_{1} <  m$ and ${\rm rank}\,{\rm {\bf B}} = r_{2} < p$, then for the
least squares solution ${\rm {\bf X}}_{LS}=(x_{ij})\in
{\mathbb{C}}^{n\times p}$  of (\ref{eq:AXB=D}) we have
\begin{equation}\label{eq:d^B}
x_{ij} = {\frac{{{\sum\limits_{\beta \in J_{r_{1},\,n} {\left\{
{i} \right\}}} { \left| {\left( {{\rm {\bf A}}^{ *} {\rm {\bf A}}}
\right)_{\,.\,i} \left( {{{\rm {\bf d}}}\,_{.\,j}^{{\rm {\bf B}}}}
\right)\, _{\beta} ^{\beta}} \right| } } }}{{{\sum\limits_{\beta
\in J_{r_{1},n}} {{\left| {\left( {{\rm {\bf A}}^{ *} {\rm {\bf
A}}} \right)_{\beta} ^{\beta} } \right|}} \sum\limits_{\alpha \in
I_{r_{2},p}}{{\left| {\left( {{\rm {\bf B}}{\rm {\bf B}}^{ *} }
\right) _{\alpha} ^{\alpha} } \right|}}} }}},
\end{equation}
or
\begin{equation}\label{eq:d^A}
 x_{ij}={\frac{{{\sum\limits_{\alpha
\in I_{r_{2},p} {\left\{ {j} \right\}}} { \left| {\left( {{\rm
{\bf B}}{\rm {\bf B}}^{ *} } \right)_{\,j\,.} \left( {{{\rm {\bf
d}}}\,_{i\,.}^{{\rm {\bf A}}}} \right)\,_{\alpha} ^{\alpha}}
\right| } }}}{{{\sum\limits_{\beta \in J_{r_{1},n}} {{\left|
{\left( {{\rm {\bf A}}^{ *} {\rm {\bf A}}} \right) _{\beta}
^{\beta} } \right|}}\sum\limits_{\alpha \in I_{r_{2},p}} {{\left|
{\left( {{\rm {\bf B}}{\rm {\bf B}}^{ *} } \right) _{\alpha}
^{\alpha} } \right|}}} }}},
\end{equation}
where
  \begin{equation}
\label{eq:def_d^B_m}
   {{\rm {\bf d}}_{.\,j}^{{\rm {\bf B}}}}=\left[
\sum\limits_{\alpha \in I_{r_{2},p} {\left\{ {j} \right\}}} {
\left| {\left( {{\rm {\bf B}}{\rm {\bf B}}^{ *} } \right)_{j.}
\left( {\tilde{{\rm {\bf d}}}_{1.}} \right)\,_{\alpha} ^{\alpha}}
\right|},...,\sum\limits_{\alpha \in I_{r_{2},p} {\left\{ {j}
\right\}}} { \left| {\left( {{\rm {\bf B}}{\rm {\bf B}}^{ *} }
\right)_{j.} \left( {\tilde{{\rm {\bf d}}}_{n.}}
\right)\,_{\alpha} ^{\alpha}} \right|} \right]^{T},
\end{equation}
  \begin{equation}
\label{eq:def_d^A_m}
  {{\rm {\bf d}}_{i\,.}^{{\rm {\bf A}}}}=\left[
\sum\limits_{\beta \in J_{r_{1},n} {\left\{ {i} \right\}}} {
\left| {\left( {{\rm {\bf A}}^{ *}{\rm {\bf A}} } \right)_{.i}
\left( {\tilde{{\rm {\bf d}}}_{.1}} \right)\,_{\beta} ^{\beta}}
\right|},...,\sum\limits_{\alpha \in I_{r_{1},n} {\left\{ {i}
\right\}}} { \left| {\left( {{\rm {\bf A}}^{ *}{\rm {\bf A}} }
\right)_{.i} \left( {\tilde{{\rm {\bf d}}}_{.\,p}}
\right)\,_{\beta} ^{\beta}} \right|} \right]
\end{equation}
 are the column-vector and the row-vector, respectively.
\item[(iii)] If ${\rm rank}\,{\rm {\bf A}} = n$ and ${\rm rank}\,{\rm {\bf B}} = r_{2} < p$, then for the
least squares solution ${\rm {\bf X}}_{LS}=(x_{ij})\in
{\mathbb{C}}^{n\times p}$  of (\ref{eq:AXB=D}) we have
\begin{equation}
\label{eq:AXB_detA*A_d^B} x_{ij}={\frac{{\det \left( {({\rm {\bf
A}}^{ *} {\rm {\bf A}})_{.\,i\,} \left( {{\rm {\bf d}}_{.j}^{{\rm
{\bf B}}}} \right)}\right)}}{{\det ({\rm {\bf A}}^{ *} {\rm {\bf
A}})\sum\limits_{\alpha \in I_{r_{2},p}}{{\left| {\left( {{\rm
{\bf B}}{\rm {\bf B}}^{ *} } \right) _{\alpha} ^{\alpha} }
\right|}}} }},
\end{equation}
or
\begin{equation}\label{AXB_detA*A_d^A}
{\frac{{{\sum\limits_{\alpha \in I_{r_{2},p} {\left\{ {j}
\right\}}} { \left| {\left( {{\rm {\bf B}}{\rm {\bf B}}^{ *} }
\right)_{\,j\,.} \left( {{{\rm {\bf d}}}\,_{i\,.}^{{\rm {\bf A}}}}
\right)\,_{\alpha} ^{\alpha}} \right| } }}}{{\det ({\rm {\bf A}}^{
*} {\rm {\bf A}})\sum\limits_{\alpha \in I_{r_{2},p}}{{\left|
{\left( {{\rm {\bf B}}{\rm {\bf B}}^{ *} } \right) _{\alpha}
^{\alpha} } \right|}}} }},
\end{equation}
where $ {{\rm {\bf d}}_{.\,j}^{{\rm {\bf B}}}}$ is
(\ref{eq:def_d^B_m})
   and
 ${\rm {\bf
d}}_{i\,.}^{{\rm {\bf A}}} $ is (\ref{eq:def_d^A}).
\item[(iiii)]If ${\rm rank}\,{\rm {\bf A}} = r_{1} <  m$ and
${\rm rank}\,{\rm {\bf B}} =  p$, then for the least squares
solution ${\rm {\bf X}}_{LS}=(x_{ij})\in {\mathbb{C}}^{n\times p}$
of (\ref{eq:AXB=D}) we have
\begin{equation}
\label{eq:AXB_detBB*_d^A} x_{i\,j} = {\frac{{\det\left( { ({\rm
{\bf B}}{\rm {\bf B}}^{ *} )_{j.\,} \left( {{\rm {\bf
d}}_{i\,.}^{{\rm {\bf A}}}} \right)}\right)}}{{{\sum\limits_{\beta
\in J_{r_{1},n}} {{\left| {\left( {{\rm {\bf A}}^{ *} {\rm {\bf
A}}} \right) _{\beta} ^{\beta} } \right|}}\cdot \det ({\rm {\bf
B}}{\rm {\bf B}}^{ *} )}}}},
\end{equation}
or
\begin{equation} \label{eq:AXB_detBB*_d^B} x_{i\,j}=
{\frac{{{\sum\limits_{\beta \in J_{r_{1},\,n} {\left\{ {i}
\right\}}} { \left| {\left( {{\rm {\bf A}}^{ *} {\rm {\bf A}}}
\right)_{\,.\,i} \left( {{{\rm {\bf d}}}\,_{.\,j}^{{\rm {\bf B}}}}
\right)\, _{\beta} ^{\beta}} \right| } } }}{{{\sum\limits_{\beta
\in J_{r_{1},n}} {{\left| {\left( {{\rm {\bf A}}^{ *} {\rm {\bf
A}}} \right)_{\beta} ^{\beta} } \right|}}\det ({\rm {\bf B}}{\rm
{\bf B}}^{ *} )}}}},
\end{equation}
 \noindent where $ {{\rm {\bf d}}_{.\,j}^{{\rm {\bf B}}}}$ is
(\ref{eq:def_d^B})
   and
 ${\rm {\bf
d}}_{i\,.}^{{\rm {\bf A}}} $ is (\ref{eq:def_d^A_m}).
\end{enumerate}
\end{theorem}
{\textit{Proof.}} i) If ${\rm rank}\,{\rm {\bf A}} = n$ and ${\rm
rank}\,{\rm {\bf B}} = p$,
  then by Corollary \ref{cor:A+A-1} ${\rm {\bf
A}}^{ +}  = \left( {{\rm {\bf A}}^{ *} {\rm {\bf A}}} \right)^{ -
1}{\rm {\bf A}}^{ *} $ and ${\rm {\bf B}}^{ +}  = {\rm {\bf B}}^{
* }\left( {{\rm {\bf B}}{\rm {\bf B}}^{ *} } \right)^{ - 1}$.
 Therefore, we obtain
\[\begin{array}{c}
{\rm {\bf X}}_{LS} = ({\rm {\bf A}}^{ \ast}{\rm {\bf A}})^{ -
1}{\rm {\bf A}}^{ \ast}{\rm {\bf D}}{\rm {\bf B}}^{ *} \left(
{{\rm {\bf
B}}{\rm {\bf B}}^{ *} } \right)^{ - 1}=\\
 = \begin{pmatrix}
   x_{11} & x_{12} & \ldots & x_{1p} \\
   x_{21} & x_{22} & \ldots & x_{2p} \\
   \ldots & \ldots & \ldots & \ldots \\
   x_{n1} & x_{n2} & \ldots & x_{np} \
 \end{pmatrix}
 = {\frac{{1}}{{\det ({\rm {\bf A}}^{ \ast}{\rm {\bf A}})}}}\begin{pmatrix}
  {L} _{11}^{{\rm {\bf A}}} & {L} _{21}^{{\rm {\bf A}}}&
   \ldots & {L} _{n1}^{{\rm {\bf A}}} \\
  {L} _{12}^{{\rm {\bf A}}} & {L} _{22}^{{\rm {\bf A}}} &
  \ldots & {L} _{n2}^{{\rm {\bf A}}} \\
  \ldots & \ldots & \ldots & \ldots \\
 {L} _{1n}^{{\rm {\bf A}}} & {L} _{2n}^{{\rm {\bf A}}}
  & \ldots & {L} _{nn}^{{\rm {\bf A}}}
\end{pmatrix}\times\\
  \times\begin{pmatrix}
    \tilde{d}_{11} & \tilde{d}_{12} & \ldots & \tilde{d}_{1m} \\
    \tilde{d}_{21} & \tilde{d}_{22} & \ldots & \tilde{d}_{2m} \\
    \ldots & \ldots & \ldots & \ldots \\
    \tilde{d}_{n1} & \tilde{d}_{n2} & \ldots & \tilde{d}_{nm} \
  \end{pmatrix}

{\frac{{1}}{{\det \left( {{\rm {\bf B}}{\rm {\bf B}}^{ *} }
\right)}}}
\begin{pmatrix}
 {R} _{\, 11}^{{\rm {\bf B}}} & {R} _{\, 21}^{{\rm {\bf B}}}
 &\ldots & {R} _{\, p1}^{{\rm {\bf B}}} \\
 {R} _{\, 12}^{{\rm {\bf B}}} & {R} _{\, 22}^{{\rm {\bf B}}} &\ldots &
 {R} _{\, p2}^{{\rm {\bf B}}} \\
 \ldots  & \ldots & \ldots & \ldots \\
 {R} _{\, 1p}^{{\rm {\bf B}}} & {R} _{\, 2p}^{{\rm {\bf B}}} &
 \ldots & {R} _{\, pp}^{{\rm {\bf B}}}
\end{pmatrix},
\end{array}
\]
\noindent where $\tilde{d}_{ij}$ is  $ij$th entry of the matrix
${\rm {\bf \widetilde{D}}}$, ${L}_{ij}^{{\rm {\bf A}}} $ is the
$ij$th
 cofactor of $({\rm {\bf A}}^{ \ast}{\rm {\bf A}})$ for all $i,j =  \overline {1,n} $ and ${R}_{i{\kern 1pt} j}^{{\rm {\bf
B}}} $ is the  $ij$th cofactor of $\left( {{\rm {\bf B}}{\rm {\bf
B}}^{ *} } \right)$ for all $i,j = \overline {1,p} $. This implies
\begin{equation}
\label{eq:sum}x_{ij} = {\frac{{{\sum\limits_{k = 1}^{n} {{
L}_{ki}^{{\rm {\bf A}}}} }\left( {{\sum\limits_{s = 1}^{p}
{\tilde{d}_{\,ks}} } {R}_{js}^{{\rm {\bf B}}}} \right)}}{{\det
({\rm {\bf A}}^{ *} {\rm {\bf A}})\cdot\det ({\rm {\bf B}}{\rm
{\bf B}}^{ *} )}}},
\end{equation}
 \noindent for all $i = \overline {1,n} $, $j =
\overline {1,p} $. We obtain the sum in parentheses  and denote it
as follows
\[{\sum\limits_{s
= 1}^{p} {\tilde{d}_{k\,s}} } {R}_{j\,s}^{{\rm {\bf B}}} = \det
({\rm {\bf B}}{\rm {\bf B}}^{ *} )_{j.\,} \left( {\tilde{{\rm {\bf
d}}}_{k\,.}} \right):=d_{k\,j}^{{\rm {\bf B}}},\] where
$\tilde{{\rm {\bf d}}}_{k\,.} $ is the $k$th row-vector of
$\tilde{{\rm {\bf D}}}$ for all $k = \overline {1,n} $. Suppose
${\rm {\bf d}}_{.\,j}^{{\rm {\bf B}}} : = \left(d_{1\,j}^{{\rm
{\bf B}}},\ldots ,d_{n\,j}^{{\rm {\bf B}}} \right)^{T}$
 is the column-vector for all
$j = \overline {1,p} $. Reducing the sum ${\sum\limits_{k = 1}^{n}
{{ L}_{ki}^{{\rm {\bf A}}}} }d_{k\,j}^{{\rm {\bf B}}} $, we obtain
an analog of Cramer's rule for (\ref{eq:AXB=D}) by
(\ref{eq:AXB_cdetA*A}).

Interchanging the order of summation in (\ref{eq:sum}), we have
\[
 x_{ij} = {\frac{{{\sum\limits_{s = 1}^{p} {\left(
{{\sum\limits_{k = 1}^{n} {{L}_{ki}^{{\rm {\bf A}}}
\tilde{d}_{\,ks}} } } \right)}} {R}_{js}^{{\rm {\bf B}}}} }{{\det
({\rm {\bf A}}^{ *} {\rm {\bf A}})\cdot\det ({\rm {\bf B}}{\rm
{\bf B}}^{ *} )}}}.\] We obtain the sum in parentheses  and denote
it as follows
\[{\sum\limits_{k = 1}^{n}
{{L}_{ki}^{{\rm {\bf A}}} \tilde{d}_{k\,s}} }  = \det ({\rm {\bf
A}}^{ *} {\rm {\bf A}})_{.\,i\,} \left( {\tilde{{\rm {\bf
d}}}_{.\,s}} \right)=:d_{i\,s}^{{\rm {\bf A}}},\] \noindent where
$\tilde{{\rm {\bf d}}}_{.\,s} $ is the $s$th column-vector of
$\tilde{{\rm {\bf D}}}$ for all $s = \overline {1,p} $. Suppose
${{{\rm {\bf d}}}_{i\,.}^{{\rm {\bf A}}}} : = \left(
d_{i\,1}^{{\rm {\bf A}}},\ldots , d_{i\,p}^{{\rm {\bf A}}}
\right)$ is the row-vector for all $i = \overline {1,n} $.
Reducing the sum ${{\sum\limits_{s = 1}^{n} d_{i\,s}^{{\rm {\bf
A}}}} {R}_{js}^{{\rm {\bf B}}}} $, we obtain another analog of
Cramer's rule for the least squares solutions of (\ref{eq:AXB=D})
by (\ref{eq:AXB_rdetBB*}).

ii) If ${\rm {\bf A}} \in {\rm {\mathbb{C}}}_{r_{1}}^{m\times n}
$, ${\rm {\bf B}} \in {\rm {\mathbb{C}}}_{r_{2}}^{p\times q} $ and
$ r_{1} <  n$, $ r_{2} < p$, then by Theorem
\ref{theor:det_repr_MP} the Moore-Penrose inverses ${\rm {\bf
A}}^{ +} = \left( {a_{ij}^{ +} } \right) \in {\rm
{\mathbb{C}}}_{}^{n\times m} $ and ${\rm {\bf B}}^{ +} = \left(
{b_{ij}^{ +} } \right) \in {\rm {\mathbb{C}}}^{q\times p} $
possess the following determinantal representations respectively,
\[
 a_{ij}^{ +}  = {\frac{{{\sum\limits_{\beta
\in J_{r_{1},\,n} {\left\{ {i} \right\}}} { \left| {\left( {{\rm
{\bf A}}^{ *} {\rm {\bf A}}} \right)_{\,. \,i} \left( {{\rm {\bf
a}}_{.j}^{ *} }  \right){\kern 1pt}  _{\beta} ^{\beta}} \right| }
} }}{{{\sum\limits_{\beta \in J_{r_{1},\,n}} {{\left| {\left(
{{\rm {\bf A}}^{ *} {\rm {\bf A}}} \right){\kern 1pt} _{\beta}
^{\beta} }  \right|}}} }}},
\]
\begin{equation}\label{eq:b+}
 b_{ij}^{ +}  =
{\frac{{{\sum\limits_{\alpha \in I_{r_{2},p} {\left\{ {j}
\right\}}} { \left| {({\rm {\bf B}}{\rm {\bf B}}^{ *} )_{j\,.\,}
({\rm {\bf b}}_{i.\,}^{ *} )\,_{\alpha} ^{\alpha}} \right| }
}}}{{{\sum\limits_{\alpha \in I_{r_{2},p}}  {{\left| {\left( {{\rm
{\bf B}}{\rm {\bf B}}^{ *} } \right){\kern 1pt} _{\alpha}
^{\alpha} } \right|}}} }}}.
\end{equation}
Since by Theorem \ref{theor:LS_AXB}  ${\rm {\bf X}}_{LS}={\rm {\bf
A}}^{+}{\rm {\bf D}}{\rm {\bf B}}^{+}$, then an entry of ${\rm
{\bf X}}_{LS}=(x_{ij})$ is

\begin{equation}
\label{eq:sum+} x_{ij} = {{\sum\limits_{s = 1}^{q} {\left(
{{\sum\limits_{k = 1}^{m} {{a}_{ik}^{+} d_{ks}} } } \right)}}
{b}_{sj}^{+}}.
\end{equation}
Denote  by $\hat{{\rm {\bf d}}_{.s}}$ the $s$th column of ${\rm
{\bf A}}^{ \ast}{\rm {\bf D}}=:\hat{{\rm {\bf D}}}=
(\hat{d}_{ij})\in {\mathbb{C}}^{n\times q}$ for all $s=\overline
{1,q}$. It follows from ${\sum\limits_{k} { {\rm {\bf a}}_{.\,k}^{
*}}d_{ks} }=\hat{{\rm {\bf d}}_{.\,s}}$ that
\[
\sum\limits_{k = 1}^{m} {{a}_{ik}^{+} d_{ks}}=\sum\limits_{k =
1}^{m}{\frac{{{\sum\limits_{\beta \in J_{r_{1},\,n} {\left\{ {i}
\right\}}} { \left| {\left( {{\rm {\bf A}}^{ *} {\rm {\bf A}}}
\right)_{\,. \,i} \left( {{\rm {\bf a}}_{.k}^{ *} } \right) {\kern
1pt} _{\beta} ^{\beta}} \right| } } }}{{{\sum\limits_{\beta \in
J_{r_{1},\,n}} {{\left| {\left( {{\rm {\bf A}}^{ *} {\rm {\bf A}}}
\right){\kern 1pt} _{\beta} ^{\beta} }  \right|}}} }}}\cdot
d_{ks}=
\]
\begin{equation}\label{eq:sum_cdet}
{\frac{{{\sum\limits_{\beta \in J_{r_{1},\,n} {\left\{ {i}
\right\}}}\sum\limits_{k = 1}^{m} { \left| {\left( {{\rm {\bf
A}}^{ *} {\rm {\bf A}}} \right)_{\,. \,i} \left( {{\rm {\bf
a}}_{.k}^{ *} } \right) {\kern 1pt} _{\beta} ^{\beta}} \right| } }
}\cdot d_{ks}}{{{\sum\limits_{\beta \in J_{r_{1},\,n}} {{\left|
{\left( {{\rm {\bf A}}^{ *} {\rm {\bf A}}} \right){\kern 1pt}
_{\beta} ^{\beta} }  \right|}}} }}}={\frac{{{\sum\limits_{\beta
\in J_{r_{1},\,n} {\left\{ {i} \right\}}} { \left| {\left( {{\rm
{\bf A}}^{ *} {\rm {\bf A}}} \right)_{\,. \,i} \left( \hat{{\rm
{\bf d}}_{.\,s}} \right) {\kern 1pt} _{\beta} ^{\beta}} \right| }
} }}{{{\sum\limits_{\beta \in J_{r_{1},\,n}} {{\left| {\left(
{{\rm {\bf A}}^{ *} {\rm {\bf A}}} \right){\kern 1pt} _{\beta}
^{\beta} }  \right|}}} }}}
\end{equation}
Suppose ${\rm {\bf e}}_{s.}$ and ${\rm {\bf e}}_{.\,s}$ are
respectively the unit row-vector and the unit column-vector whose
components are $0$, except the $s$th components, which are $1$.
Substituting  (\ref{eq:sum_cdet}) and (\ref{eq:b+}) in
(\ref{eq:sum+}), we obtain
\[
x_{ij} =\sum\limits_{s = 1}^{q}{\frac{{{\sum\limits_{\beta \in
J_{r_{1},\,n} {\left\{ {i} \right\}}} { \left| {\left( {{\rm {\bf
A}}^{ *} {\rm {\bf A}}} \right)_{\,. \,i} \left( \hat{{\rm {\bf
d}}_{.\,s}} \right) {\kern 1pt} _{\beta} ^{\beta}} \right| } }
}}{{{\sum\limits_{\beta \in J_{r_{1},\,n}} {{\left| {\left( {{\rm
{\bf A}}^{ *} {\rm {\bf A}}} \right){\kern 1pt} _{\beta} ^{\beta}
}  \right|}}} }}}{\frac{{{\sum\limits_{\alpha \in I_{r_{2},p}
{\left\{ {j} \right\}}} { \left| {({\rm {\bf B}}{\rm {\bf B}}^{ *}
)_{j\,.\,} ({\rm {\bf b}}_{s.\,}^{ *} )\,_{\alpha} ^{\alpha}}
\right| } }}}{{{\sum\limits_{\alpha \in I_{r_{2},p}}  {{\left|
{\left( {{\rm {\bf B}}{\rm {\bf B}}^{ *} } \right){\kern 1pt}
_{\alpha} ^{\alpha} } \right|}}} }}}.
\]
Since
\begin{equation}\label{eq:prop}
\hat{{\rm {\bf d}}_{.\,s}}=\sum\limits_{l = 1}^{n}{\rm {\bf
e}}_{.\,l}\hat{ d_{ls}},\,\,\,\,  {\rm {\bf b}}_{s.\,}^{
*}=\sum\limits_{t = 1}^{p}b_{st}^{*}{\rm {\bf e}}_{t.},\,\,\,\,
\sum\limits_{s=1}^{q}\hat{d_{ls}}b_{st}^{*}=\widetilde{d}_{lt},
\end{equation}
then we have
\[
x_{ij} = \]
\[{\frac{{ \sum\limits_{s = 1}^{q}\sum\limits_{t =
1}^{p} \sum\limits_{l = 1}^{n} {\sum\limits_{\beta \in
J_{r_{1},\,n} {\left\{ {i} \right\}}} { \left| {\left( {{\rm {\bf
A}}^{ *} {\rm {\bf A}}} \right)_{\,. \,i} \left( {\rm {\bf
e}}_{.\,l} \right) {\kern 1pt} _{\beta} ^{\beta}} \right| } }
}\hat{ d_{ls}}b_{st}^{*}{\sum\limits_{\alpha \in I_{r_{2},p}
{\left\{ {j} \right\}}} { \left| {({\rm {\bf B}}{\rm {\bf B}}^{ *}
)_{j\,.\,} ({\rm {\bf e}}_{t.} )\,_{\alpha} ^{\alpha}} \right| } }
}{{{\sum\limits_{\beta \in J_{r_{1},\,n}} {{\left| {\left( {{\rm
{\bf A}}^{ *} {\rm {\bf A}}} \right){\kern 1pt} _{\beta} ^{\beta}
}  \right|}}} }{{\sum\limits_{\alpha \in I_{r_{2},p}} {{\left|
{\left( {{\rm {\bf B}}{\rm {\bf B}}^{ *} } \right){\kern 1pt}
_{\alpha} ^{\alpha} } \right|}}} }}    }=
\]
\begin{equation}\label{eq:x_ij}
{\frac{{ \sum\limits_{t = 1}^{p} \sum\limits_{l = 1}^{n}
{\sum\limits_{\beta \in J_{r_{1},\,n} {\left\{ {i} \right\}}}
{\left| {\left( {{\rm {\bf A}}^{ *} {\rm {\bf A}}} \right)_{\,.
\,i} \left( {\rm {\bf e}}_{.\,l} \right) {\kern 1pt} _{\beta}
^{\beta}} \right| } } }\,\,\widetilde{d}_{lt}{\sum\limits_{\alpha
\in I_{r_{2},p} {\left\{ {j} \right\}}} { \left| {({\rm {\bf
B}}{\rm {\bf B}}^{ *} )_{j\,.\,} ({\rm {\bf e}}_{t.} )\,_{\alpha}
^{\alpha}} \right| } } }{{{\sum\limits_{\beta \in J_{r_{1},\,n}}
{{\left| {\left( {{\rm {\bf A}}^{ *} {\rm {\bf A}}} \right){\kern
1pt} _{\beta} ^{\beta} }  \right|}}} }{{\sum\limits_{\alpha \in
I_{r_{2},p}} {{\left| {\left( {{\rm {\bf B}}{\rm {\bf B}}^{ *} }
\right){\kern 1pt} _{\alpha} ^{\alpha} } \right|}}} }}    }.
\end{equation}
Denote by
\[
 d^{{\rm {\bf A}}}_{it}:= \]
\[
{\sum\limits_{\beta \in J_{r_{1},\,n} {\left\{ {i} \right\}}} {
\left| {\left( {{\rm {\bf A}}^{ *} {\rm {\bf A}}} \right)_{\,.
\,i} \left( \widetilde{{\rm {\bf d}}}_{.\,t} \right){\kern 1pt}
_{\beta} ^{\beta}} \right| } }= \sum\limits_{l = 1}^{n}
{\sum\limits_{\beta \in J_{r_{1},\,n} {\left\{ {i} \right\}}} {
\left| {\left( {{\rm {\bf A}}^{ *} {\rm {\bf A}}} \right)_{\,.
\,i} \left( {\rm {\bf e}}_{.\,l} \right){\kern 1pt}
 _{\beta} ^{\beta}} \right| } } \widetilde{d}_{lt}
\]
the $t$th component  of a row-vector ${\rm {\bf d}}^{{\rm {\bf
A}}}_{i\,.}= (d^{{\rm {\bf A}}}_{i1},...,d^{{\rm {\bf A}}}_{ip})$
for all $t=\overline {1,p}$. Substituting it in (\ref{eq:x_ij}),
we have
\[x_{ij} ={\frac{{ \sum\limits_{t = 1}^{p}
 d^{{\rm {\bf A}}}_{it}
}{\sum\limits_{\alpha \in I_{r_{2},p} {\left\{ {j} \right\}}} {
\left| {({\rm {\bf B}}{\rm {\bf B}}^{ *} )_{j\,.\,} ({\rm {\bf
e}}_{t.} )\,_{\alpha} ^{\alpha}} \right| } }
}{{{\sum\limits_{\beta \in J_{r_{1},\,n}} {{\left| {\left( {{\rm
{\bf A}}^{ *} {\rm {\bf A}}} \right){\kern 1pt} _{\beta} ^{\beta}
}  \right|}}} }{{\sum\limits_{\alpha \in I_{r_{2},p}} {{\left|
{\left( {{\rm {\bf B}}{\rm {\bf B}}^{ *} } \right){\kern 1pt}
_{\alpha} ^{\alpha} } \right|}}} }}    }.
\]
Since $\sum\limits_{t = 1}^{p}
 d^{{\rm {\bf A}}}_{it}{\rm {\bf e}}_{t.}={\rm {\bf
d}}^{{\rm {\bf A}}}_{i\,.}$, then it follows (\ref{eq:d^A}).

If we denote by
\begin{equation}\label{eq:d^B_den}
 d^{{\rm {\bf B}}}_{lj}:=
\sum\limits_{t = 1}^{p}\widetilde{d}_{lt}{\sum\limits_{\alpha \in
I_{r_{2},p} {\left\{ {j} \right\}}} { \left| {({\rm {\bf B}}{\rm
{\bf B}}^{ *} )_{j\,.\,} ({\rm {\bf e}}_{t.} )\,_{\alpha}
^{\alpha}} \right| } }={\sum\limits_{\alpha \in I_{r_{2},p}
{\left\{ {j} \right\}}} {\left| {({\rm {\bf B}}{\rm {\bf B}}^{ *}
)_{j\,.\,} (\widetilde{{\rm {\bf d}}}_{l.} )\,_{\alpha} ^{\alpha}}
\right| } }
\end{equation}
the $l$th component  of a column-vector ${\rm {\bf d}}^{{\rm {\bf
B}}}_{.\,j}= (d^{{\rm {\bf B}}}_{1j},...,d^{{\rm {\bf
B}}}_{jn})^{T}$ for all $l=\overline {1,n}$ and substitute it in
(\ref{eq:x_ij}), we obtain
\[x_{ij} ={\frac{{  \sum\limits_{l = 1}^{n}
{\sum\limits_{\beta \in J_{r_{1},\,n} {\left\{ {i} \right\}}} {
\left| {\left( {{\rm {\bf A}}^{ *} {\rm {\bf A}}} \right)_{\,.
\,i} \left( {\rm {\bf e}}_{.\,l} \right){\kern 1pt} _{\beta}
^{\beta}} \right| } } }\,\,d^{{\rm {\bf B}}}_{lj}
}{{{\sum\limits_{\beta \in J_{r_{1},\,n}} {{\left| {\left( {{\rm
{\bf A}}^{ *} {\rm {\bf A}}} \right){\kern 1pt} _{\beta} ^{\beta}
}  \right|}}} }{{\sum\limits_{\alpha \in I_{r_{2},p}} {{\left|
{\left( {{\rm {\bf B}}{\rm {\bf B}}^{ *} } \right){\kern 1pt}
_{\alpha} ^{\alpha} } \right|}}} }}    }.
\]
Since $\sum\limits_{l = 1}^{n}{\rm {\bf e}}_{.l}
 d^{{\rm {\bf B}}}_{lj}={\rm {\bf
d}}^{{\rm {\bf B}}}_{.\,j}$, then it follows (\ref{eq:d^B}).

iii) If ${\rm {\bf A}} \in {\rm {\mathbb{C}}}_{r_{1}}^{m\times n}
$, ${\rm {\bf B}} \in {\rm {\mathbb{C}}}_{r_{2}}^{p\times q} $ and
$ r_{1} =n$, $ r_{2} < p$, then by Theorem \ref{theor:det_repr_MP}
and  Remark \ref{rem:det_repr_MP_A*A} the Moore-Penrose inverses
${\rm {\bf A}}^{ +} = \left( {a_{ij}^{ +} } \right) \in {\rm
{\mathbb{C}}}_{}^{n\times m} $ and ${\rm {\bf B}}^{ +} = \left(
{b_{ij}^{ +} } \right) \in {\rm {\mathbb{C}}}^{q\times p} $
possess the following determinantal representations respectively,
\[
 a_{ij}^{ +}  = {\frac{\det {\left( {{\rm
{\bf A}}^{ *} {\rm {\bf A}}} \right)_{\,. \,i} \left( {{\rm {\bf
a}}_{.j}^{ *} }  \right)}   } {\det {\left( {{\rm {\bf A}}^{ *}
{\rm {\bf A}}} \right)}}},
\]
\begin{equation}\label{eq:b+2}
 b_{ij}^{ +}  =
{\frac{{{\sum\limits_{\alpha \in I_{r_{2},p} {\left\{ {j}
\right\}}} { \left| {({\rm {\bf B}}{\rm {\bf B}}^{ *} )_{j\,.\,}
({\rm {\bf b}}_{i.\,}^{ *} )\,_{\alpha} ^{\alpha}} \right| }
}}}{{{\sum\limits_{\alpha \in I_{r_{2},p}}  {{\left| {\left( {{\rm
{\bf B}}{\rm {\bf B}}^{ *} } \right){\kern 1pt} _{\alpha}
^{\alpha} } \right|}}} }}}.
\end{equation}
Since by Theorem \ref{theor:LS_AXB}  ${\rm {\bf X}}_{LS}={\rm {\bf
A}}^{+}{\rm {\bf D}}{\rm {\bf B}}^{+}$, then an entry of ${\rm
{\bf X}}_{LS}=(x_{ij})$ is (\ref{eq:sum+}). Denote  by $\hat{{\rm
{\bf d}}_{.s}}$ the $s$th column of ${\rm {\bf A}}^{ \ast}{\rm
{\bf D}}=:\hat{{\rm {\bf D}}}= (\hat{d}_{ij})\in
{\mathbb{C}}^{n\times q}$ for all $s=\overline {1,q}$. It follows
from ${\sum\limits_{k} { {\rm {\bf a}}_{.\,k}^{ *}}d_{ks}
}=\hat{{\rm {\bf d}}_{.\,s}}$ that
\begin{equation}\label{eq:sum_det}
\sum\limits_{k = 1}^{m} {{a}_{ik}^{+} d_{ks}}=\sum\limits_{k =
1}^{m}{\frac{\det {\left( {{\rm {\bf A}}^{ *} {\rm {\bf A}}}
\right)_{\,. \,i} \left( {{\rm {\bf a}}_{.k}^{ *} }  \right)}   }
{\det {\left( {{\rm {\bf A}}^{ *} {\rm {\bf A}}} \right)}}}\cdot
d_{ks}={\frac{\det {\left( {{\rm {\bf A}}^{ *} {\rm {\bf A}}}
\right)_{\,. \,i} \left( {\hat{{\rm {\bf d}}_{.\,s}} }  \right)} }
{\det {\left( {{\rm {\bf A}}^{ *} {\rm {\bf A}}} \right)}}}
\end{equation}
Substituting  (\ref{eq:sum_det}) and (\ref{eq:b+2}) in
(\ref{eq:sum+}), and using (\ref{eq:prop}) we have
\[
x_{ij} =\sum\limits_{s = 1}^{q}{\frac{\det {\left( {{\rm {\bf
A}}^{ *} {\rm {\bf A}}} \right)_{\,. \,i} \left( {\hat{{\rm {\bf
d}}_{.\,s}} }  \right)} } {\det {\left( {{\rm {\bf A}}^{ *} {\rm
{\bf A}}} \right)}}}{\frac{{{\sum\limits_{\alpha \in I_{r_{2},p}
{\left\{ {j} \right\}}} { \left| {({\rm {\bf B}}{\rm {\bf B}}^{ *}
)_{j\,.\,} ({\rm {\bf b}}_{s.\,}^{ *} )\,_{\alpha} ^{\alpha}}
\right| } }}}{{{\sum\limits_{\alpha \in I_{r_{2},p}} {{\left|
{\left( {{\rm {\bf B}}{\rm {\bf B}}^{ *} } \right){\kern 1pt}
_{\alpha} ^{\alpha} } \right|}}} }}}=
\]
\[{\frac{{ \sum\limits_{s = 1}^{q}\sum\limits_{t =
1}^{p} \sum\limits_{l = 1}^{n} \det {\left( {{\rm {\bf A}}^{ *}
{\rm {\bf A}}} \right)_{\,. \,i} \left( {{\rm {\bf e}}_{.\,l} }
\right)} }\hat{ d_{ls}}b_{st}^{*}{\sum\limits_{\alpha \in
I_{r_{2},p} {\left\{ {j} \right\}}} { \left| {({\rm {\bf B}}{\rm
{\bf B}}^{ *} )_{j\,.\,} ({\rm {\bf e}}_{t.} )\,_{\alpha}
^{\alpha}} \right| } } }{{{\det {\left( {{\rm {\bf A}}^{ *} {\rm
{\bf A}}} \right)}} }{{\sum\limits_{\alpha \in I_{r_{2},p}}
{{\left| {\left( {{\rm {\bf B}}{\rm {\bf B}}^{ *} } \right){\kern
1pt} _{\alpha} ^{\alpha} } \right|}}} }}    }=
\]
\begin{equation}\label{eq:x_ij_2}
{\frac{{ \sum\limits_{t = 1}^{p} \sum\limits_{l = 1}^{n} \det
{\left( {{\rm {\bf A}}^{ *} {\rm {\bf A}}} \right)_{\,. \,i}
\left( {{\rm {\bf e}}_{.\,l} } \right)}
}\,\,\widetilde{d}_{lt}{\sum\limits_{\alpha \in I_{r_{2},p}
{\left\{ {j} \right\}}} { \left| {({\rm {\bf B}}{\rm {\bf B}}^{ *}
)_{j\,.\,} ({\rm {\bf e}}_{t.} )\,_{\alpha} ^{\alpha}} \right| } }
}{{{\det {\left( {{\rm {\bf A}}^{ *} {\rm {\bf A}}} \right)}}
}{{\sum\limits_{\alpha \in I_{r_{2},p}} {{\left| {\left( {{\rm
{\bf B}}{\rm {\bf B}}^{ *} } \right){\kern 1pt} _{\alpha}
^{\alpha} } \right|}}} }}    }.
\end{equation}
If we integrate (\ref{eq:d^B_den}) in (\ref{eq:x_ij_2}), then we
get
\[x_{ij} ={\frac{{  \sum\limits_{l = 1}^{n}
\det {\left( {{\rm {\bf A}}^{ *} {\rm {\bf A}}} \right)_{\,. \,i}
\left( {{\rm {\bf e}}_{.\,l} } \right)} }\,\,d^{{\rm {\bf
B}}}_{lj} }{{{\det {\left( {{\rm {\bf A}}^{ *} {\rm {\bf A}}}
\right)}} }{{\sum\limits_{\alpha \in I_{r_{2},p}} {{\left| {\left(
{{\rm {\bf B}}{\rm {\bf B}}^{ *} } \right){\kern 1pt} _{\alpha}
^{\alpha} } \right|}}} }}    }.
\]
Since again $\sum\limits_{l = 1}^{n}{\rm {\bf e}}_{.l}
 d^{{\rm {\bf B}}}_{lj}={\rm {\bf
d}}^{{\rm {\bf B}}}_{.\,j}$, then it follows
(\ref{eq:AXB_detA*A_d^B}), where $ {{\rm {\bf d}}_{.\,j}^{{\rm
{\bf B}}}}$ is (\ref{eq:def_d^B_m}).

If we denote by
\[
 d^{{\rm {\bf A}}}_{it}:= \]
\[
\sum\limits_{l = 1}^{n} \det {\left( {{\rm {\bf A}}^{ *} {\rm {\bf
A}}} \right)_{\,. \,i} \left( {\widetilde{{\rm {\bf d}}}_{.\,t} }
\right)} = \sum\limits_{l = 1}^{n} \det {\left( {{\rm {\bf A}}^{
*} {\rm {\bf A}}} \right)_{\,. \,i} \left( {{\rm {\bf e}}_{.\,l} }
\right)} \,\,\widetilde{d}_{lt}
\]
the $t$th component  of a row-vector ${\rm {\bf d}}^{{\rm {\bf
A}}}_{i\,.}= (d^{{\rm {\bf A}}}_{i1},...,d^{{\rm {\bf A}}}_{ip})$
for all $t=\overline {1,p}$ and substitute it in
(\ref{eq:x_ij_2}), we obtain
\[x_{ij} ={\frac{{ \sum\limits_{t = 1}^{p}
 d^{{\rm {\bf A}}}_{it}
}{\sum\limits_{\alpha \in I_{r_{2},p} {\left\{ {j} \right\}}} {
\left| {({\rm {\bf B}}{\rm {\bf B}}^{ *} )_{j\,.\,} ({\rm {\bf
e}}_{t.} )\,_{\alpha} ^{\alpha}} \right| } } }{{{\det {\left(
{{\rm {\bf A}}^{ *} {\rm {\bf A}}} \right)}}
}{{\sum\limits_{\alpha \in I_{r_{2},p}} {{\left| {\left( {{\rm
{\bf B}}{\rm {\bf B}}^{ *} } \right){\kern 1pt} _{\alpha}
^{\alpha} } \right|}}} }}    }.
\]
Since again $\sum\limits_{t = 1}^{p}
 d^{{\rm {\bf A}}}_{it}{\rm {\bf e}}_{t.}={\rm {\bf
d}}^{{\rm {\bf A}}}_{i\,.}$, then it follows
(\ref{AXB_detA*A_d^A}), where ${\rm {\bf d}}^{{\rm {\bf
A}}}_{i\,.}$ is (\ref{eq:def_d^A}).

iiii) The proof is similar to the proof of iii). $\blacksquare$

\section{An example}
In this section, we give an example to illustrate our results. Let
us consider the matrix equation
\begin{equation}\label{eq_ex:AXB=D}
 {\rm {\bf A}}{\rm {\bf X}}{\rm {\bf B}} = {\rm {\bf
D}},
\end{equation}
where
\[{\bf A}=\begin{pmatrix}
  1 & i & i \\
  i & -1 & -1 \\
  0 & 1 & 0 \\
  -1 & 0 & -i
\end{pmatrix},\,\, {\bf B}=\begin{pmatrix}
  i & 1 & -i \\
  -1 & i & 1
\end{pmatrix},\,\, {\bf D}=\begin{pmatrix}
  1 & i & 1 \\
  i & 0 & 1 \\
  1 & i & 0 \\
  0 & 1 & i
\end{pmatrix}.\]
Since ${\rm rank}\,{\rm {\bf A}} = 2$ and ${\rm rank}\,{\rm {\bf
B}} = 1$, then we have the case (ii) of Theorem \ref{theor:AXB=D}.
We shall find the least squares solution of (\ref{eq_ex:AXB=D}) by
(\ref{eq:d^B}). Then we have

\[{\bf A}^{*}{\bf A}=\begin{pmatrix}
  3 & 2i & 3i \\
  -2i & 3 & 2 \\
  -3i & 2 & 3
\end{pmatrix},\,\, {\bf B}{\bf B}^{*}=\begin{pmatrix}
  3 & -3i \\
  3i & 3
\end{pmatrix},\,\, {\rm {\bf \widetilde{D}}}= {\rm {\bf
A}}^\ast{\rm {\bf D}}{\rm {\bf B}}^{\ast}=\begin{pmatrix}
  1 & -i \\
  -i & -1 \\
  -i & -1
\end{pmatrix},\] and $
{{{\sum\limits_{\alpha \in I_{1,\,2}}  {{\left| {\left( {{\rm {\bf
B}}{\rm {\bf B}}^{ *} } \right) {\kern 1pt} _{\alpha} ^{\alpha} }
\right|}}} }}=3+3=6,$
\[ {{{\sum\limits_{\beta \in J_{2,\,3}}
{{\left| {\left( {{\rm {\bf A}}^{ *}{\rm {\bf A}} } \right) {\kern
1pt} _{\beta} ^{\beta} } \right|}}} }}=\det\begin{pmatrix}
  3 & 2i \\
  -2i & 3
\end{pmatrix}+\det\begin{pmatrix}
  3 & 2 \\
  2 & 3
\end{pmatrix}+\det\begin{pmatrix}
  3 & 3i \\
  -3i & 3
\end{pmatrix}=12.\]
By (\ref{eq:def_d^B}), we can get \[{\bf d}_{.1}^{{\bf
B}}=\begin{pmatrix}
  1 \\
  -i \\
  -i
\end{pmatrix},\,\,\,\,\,{\bf d}_{.2}^{{\bf  B}}=\begin{pmatrix}
  -i \\
  -1 \\
  -1
\end{pmatrix}.\]
Since $\left( {{\rm {\bf A}}^{ *} {\rm {\bf A}}} \right)_{\,.\,1}
\left( {{{\rm {\bf d}}}\,_{.\,1}^{{\rm {\bf B}}}} \right)=
\begin{pmatrix}
  1 & 2i & 3i \\
  -i & 3 & 2 \\
  -i & 2 & 3
\end{pmatrix}$, then finally we obtain
\[
x_{11} = {\frac{{{\sum\limits_{\beta \in J_{2,\,3} {\left\{ {i}
\right\}}} { \left| {\left( {{\rm {\bf A}}^{ *} {\rm {\bf A}}}
\right)_{\,.\,1} \left( {{{\rm {\bf d}}}\,_{.\,1}^{{\rm {\bf B}}}}
\right)\, _{\beta} ^{\beta}} \right| } } }}{{{\sum\limits_{\beta
\in J_{2,3}} {{\left| {\left( {{\rm {\bf A}}^{ *} {\rm {\bf A}}}
\right)_{\beta} ^{\beta} } \right|}} \sum\limits_{\alpha \in
I_{1,2}}{{\left| {\left( {{\rm {\bf B}}{\rm {\bf B}}^{ *} }
\right) _{\alpha} ^{\alpha} } \right|}}}
}}}=\frac{\det\begin{pmatrix}
  1 & 2i \\
  -i & 3
\end{pmatrix}+\det\begin{pmatrix}
  1 & 3i \\
  -i & 3
\end{pmatrix}}{72}=-\frac{1}{72}.\]
Similarly,
\[
  x_{12} =\frac{\det\begin{pmatrix}
  -i & 2i \\
  -1 & 3
\end{pmatrix}+\det\begin{pmatrix}
  -i & 3i \\
  -1 & 3
\end{pmatrix}}{72}=-\frac{i}{72},\]    \[x_{21} =\frac{\det\begin{pmatrix}
  3 & 1 \\
  -2i & -i
\end{pmatrix}+\det\begin{pmatrix}
  -i & 2 \\
  -i & 3
\end{pmatrix}}{72}=-\frac{2i}{72},\]
\[
   x_{22} =\frac{\det\begin{pmatrix}
  3 & -i \\
  -2i & -1
\end{pmatrix}+\det\begin{pmatrix}
  -1 & 2 \\
  -1 & 3
\end{pmatrix}}{72}=-\frac{2}{72},\] \[   x_{31} =\frac{\det\begin{pmatrix}
  3 & 1 \\
  -3i & -i
\end{pmatrix}+\det\begin{pmatrix}
  3 & -i \\
  2 & -i
\end{pmatrix}}{72}=-\frac{i}{72},\]
\[
    x_{32} =\frac{\det\begin{pmatrix}
  3 & -i \\
  -3i & -1
\end{pmatrix}+\det\begin{pmatrix}
  3 & -1 \\
  2 & -1
\end{pmatrix}}{72}=-\frac{1}{72}.
\]

\end{document}